\documentclass[12pt]{article}

\newcommand{\br}{{\bf R}}

\newcommand{\bp}{{\bf P}}
\newcommand{\bq}{{\bf Q}}
\newcommand{\bz}{{\bf Z}}

\newcommand{\om}{{\omega}}
\newcommand{\la}{\lambda}
\newcommand{\bs}{{\bf S}}
\newcommand{\cm}{{\cal M}}
\newcommand{\cs}{{\cal S}}

\usepackage{latexsym,amsfonts,amsmath,amsthm,amssymb,graphicx}

\usepackage[T1]{fontenc}
\usepackage[latin1]{inputenc}
\title{Generalized probabilities taking values in non-Archimedean fields and topological groups}
\author{Andrei Khrennikov\\
School of Mathematics and Systems Engineering\\
University of V\"axj\"o, S-35195, Sweden}

\begin{document}

\maketitle

\abstract{We develop an analogue of probability theory for
probabilities taking values in topological groups. We generalize
Kolmogorov's method of axiomatization of probability theory: main
distinguishing  features of frequency probabilities are taken as
axioms in the measure-theoretic approach. We also present a review
of non-Kolmogorovian probabilistic models including models with
negative, complex, and $p$-adic   valued probabilities. The latter
model is discussed in details. The introduction of $p$-adic (as well
as more general non-Archimedean) probabilities is one of the main
motivations for consideration of generalized probabilities taking
values in topological groups which are distinct from the field of
real numbers. We discuss applications of non-Kolmogorovian models in
physics and cognitive sciences. An important part of this paper is
devoted to statistical interpretation of probabilities taking values
in topological groups (and in particular in non-Archimedean
fields).}

\section{Introduction}
Since the creation of the modern probabilistic axiomatics by A. N.
Kolmogorov [1] in 1933, probability theory was reduced to the theory
of normalized $\sigma$-additive measures taking values in the
segment [0,1] of the field of real numbers $\br$. In particular, the
main competitor of Kolmogorov's measure-theoretic approach, von
Mises' frequency approach to probability [2], [3], practically
totally disappeared from the probabilistic arena. On one hand, this
was a consequence of difficulties with von Mises' definition of
randomness (via place selections), see e.g.,
[4]-[10]\footnote{However, see also [11]-[13], where von Mises'
approach was simplified, generalized, and then fruitfully applied to
theoretical physics.}. On the other hand, von Mises' approach (as
many others) could not compete with precisely and simply formulated
Kolmogorov's theory.

We mentioned von Mises' approach not only, because its attraction
for applications, but also because von Mises' model with frequency
probabilities played the important role in the process of
formulation of the conventional axiomatics of probability theory. If
one opens Kolmogorov's book [1], he will see numerous remarks about
von Mises' theory. Andrei Nikolaevich Kolmogorov used properties of
the frequency probability to justify his choice of the axioms for
probability. In particular, Kolmogorov's probability belongs to the
segment [0,1] of the real line ${\bf R}$, because the same takes
place for von Mises' frequency probability (frequencies $\nu_N=n/N$
as well as their limits always belong to  the segment [0,1] of the
real line ${\bf R})$. In the same way Kolmogorov's probability is
additive, because the frequency probability is additive: the limit
of the sum of two frequencies equals to the sum of limits. And so
on... Thus by using {\it THEOREMS} of von Mises' frequency theory
Kolmogorov justified the AXIOMATIZATION of probability as a
normalized finite-additive measure taking values in [0,1]. Finally,
he added the condition of $\sigma$-additivity. However, the latter
condition is a purely mathematical technical condition that provides
the fruitful theory of integration (to define mean values of random
variables), see Kolmogorov remark on this condition [1]. The main
lesson of Kolmogorov's axiomatization of modern probability theory
is that {\it THEOREMS} of the frequency  probability theory were
transformed into {\it AXIOMS} in the measure-theoretic probability
theory. The Kolmogorov's axiomatics is based on:

(a) the measure-theoretic formalization of properties of relative
frequencies;

(b) Lebesque integration theory.

Regarding to (a), we would like to mention that Kolmogorov's (as
well as von Mises') assumptions were also based on a fundamental,
but hidden, assumption:

\medskip

{\it Limiting behaviour of relative frequencies is considered with
respect to one fixed topology on the field of rational numbers
$\bq,$ namely the real topology.}

\medskip

In particular, the consideration of this asymptotic behaviour implies that probabilities belong to the
 field real numbers $\br.$ In fact,  additivity of the probability is a consequence of the fact that $\br$
 is an {\it additive topological group.} We also remark that Bayes' formula
 $$
 \bp(B\vert A)=\frac{\bp(A \cap B)}{\bp(A)}, \bp(A) \ne 0.
 $$
is also a theorem in von Mises' theory, see, e.g., [11]. It is
derived as a consequence of the fact that $\br \setminus \{ 0 \}$ is
a {\it multiplicative topological group.} We can ask ourself: Are we
satisfied by Kolmogorov's theory? As pure mathematician, I would
definitely answer: "yes!". However, as physicist, I would not be so
optimistic. It seems that Kolmogorov's model, despite its
generality, does not provide a reasonable mathematical description
of all probabilistic structures that appear in physics (as well as
other natural and social sciences). In particular, we can recall the
old problem of {\it negative probabilities}. There are many objects
that by their physical origin must be probabilities, but can take
negative values (as well as values larger than 1), see e.g., Dirac
[14] (quantization of electromagnetic field), Wigner [15] (phase
space distribution), Feynman [16] (computer simulation of quantum
reality), see also [17], [11], [18].

As mathematicians, we deny these ``negative probabilities'' (because
such objects do not belong to the domain of Kolmogorov's theory). As
a consequence, physicists should work with such objects on teh
physical level of rigorousness. But negative probabilities appear
again and again in different domains of physics. The same situation
takes place for {\it complex probabilities}, see e.g., Dirac [19] or
Prugovesky [20].

We also pay attention to another probability-like structure
that recently appeared in theoretical physics. This is so called {\it $p$-adic probability}.

We recall that the first $p$-adic physical model, namely $p$-adic string,
was proposed by I. Volovich [21]. His paper induced the storm of publications
on $p$-adic string theory, e.g., [22]-[24], general $p$-adic quantum physics, e.g., [24-27],
applications of $p$-adic numbers to foundations of conventional quantum physics
(Bell's inequality, Einstein-Podolsky-Rosen paradox) [28], [29], dynamical systems [26], [30]-[32],
biological and cognitive models [26], [33].

In fact, there are two types of $p$-adic physical models:

(A) variables are $p$-adic, but functions are ${\bf C}$-valued;

(B) both variables and functions take $p$-adic values.

The A-models of $p$-adic physics and their relation to conventional
probability theory on locally compact groups (especially, totally
disconnected) will be briefly discussed in appendix 1. The B-models
are the most interesting for our present considerations.

Here we have quantities that have to be probabilities (by their
physical origin), but belong to fields of $p$-adic numbers $\bq_p$.
We again could not use Kolmogorov's axiomatics, see [25]. We also
can mention {\it comparative probability } theory, see e.g., Fine
[34], as an example of non-Kolmogorovian probabilistic model. Thus
the present situation in probability theory has some similarities
with the situation in geometry in the 19th century. We have to
recognize that Kolmogorov's theory is just one of numerous
probabilistic models. Besides Kolmogorov's model, there are numerous
non-Kolmogorovian models that describe various probabilistic
phenomena, see Accardi [35]. Such non-Kolmogorovian models can be
developed in many ways.

We mention one slight (but very important in quantum physics)
modification of Kolmogorov's theory. L. Accardi [35] proposed to
exclude from the axiomatics of probability theory Bayes' formula,
Kolmogorov's definition of conditional probability. Thus everything
is as in Kolmogorov's model besides conditional probability. L.
Accardi demonstrated that such a model could be used in the
connection with Bell's inequality and Einstein-Podolsky-Rosen
paradox.

We can construct another  large class of non-Kolmogorovian models by
considering ``probabilistic measures'' that are not defined on
$\sigma$-fields see e.g., the theory of probabilistic manifolds of
Gudder [36] that was successfully applied in quantum physics. In
fact, A. Kolmogorov considered in his first variant of axiomatics of
probability theory ``probabilities'' that are not defined on a
$\sigma$-field (and even a field) [37]. In particular, so called
density of natural numbers was considered  as probability in [37].

In this paper we shall concentrate our study to probabilistic models
that could be obtained through changing the range of values of
probabilities. Thus our ``generalized probabilities'' do not more
take values in the segment [0,1] of $\br.$

There are many ways to develop such probabilistic models. One class
of such models is related to comparative probability theory, [34].
We choose another approach for modifying the Kolmogorov's
axiomatics. We obtain natural generalizations of properties of
probabilities which are induced by the transition from $\br$ to an
arbitrary topological group. We consider $\br$ as a topological
group (with respect to addition) and extract the main properties of
Kolmogorov's measure-theoretic or von Mises' frequency probability
corresponding to the group structure on $\br.$ Then we use
generalizations of such properties to define generalized
probabilities taking values in an arbitrary topological group $G.$

On one hand, such an extension of probability  can have a lot of
applications (in particular, justification of negative or $p$-adic
probabilities). On the other hand, it could revolutionize the
classical theory of {\it Probability on Topological Structures} by
generating a huge class of new purely mathematical problems.

Before developing the general axiomatics, we  will present an
extended review on  the $p$-adic valued probabilities [25], [26],
[38], [39]. In fact, $p$-adic probability theory  was the first
example of the mathematically rigorous formalism for probabilities
taking values in a topological group $G$ which is different from
$\br.$ Since in applications of probability ( e.g., to physics) its
interpretation plays the fundamental role, we pay the main attention
to a statistical interpretation of such generalized probabilities,
see section 6.

Finally, we pay attention that recently an extremely interesting
analysis of foundations of probability theory was done by Viktor
Pavlovich Maslov in the process of development of {\it ultra-second
quantization} [39] as well as applications of classical and quantum
probabilistic models and their generalizations in finances (private
communication of V. P. Maslov).

\section {$p$-adic lessons}
\subsection{$p$-adic numbers}
The field of real numbers $\br$ is constructed as the completion of
the field of rational numbers $\bq$ with respect to the metric $p(x,
y)= |x-y|,$ where $|\cdot|$ is the usual valuation given by the
absolute value. The fields of $p$-adic numbers $\bq_p$ are
constructed in a corresponding way, but by using other valuations.
For a prime number $p$ the $p$-adic valuation $|\cdot|_p$ is defined
in the following way. First we define it for natural numbers. Every
natural number $n$ can be represented as the product of prime
numbers, $n=2^{r_2}3^{r_3}\ldots p^{r_p}\ldots,$ and we define
$|n|_p=p^{-r_p},$ writing $|0|_p=0$ and $|-n|_p=|n|_p.$ We then
extend the definition of the $p$-adic valuation $|\cdot|_p$ to all
rational numbers by setting $|n/m|_p=|n|_p/\vert m|_p$ for $m \ne
0.$ The completion of $\bq$ with respect to the metric $\rho_p(x,
y)=|x-y|_p$ is the locally compact field of $p$-adic numbers
$\bq_p.$

The number fields $\br$ and $\bq_p$ are unique in a  sense, since by
Ostrovsky's theorem, see e.g., [40], $|\cdot|$ and $|\cdot|_p$ are
the only possible valuations on $\bq,$ but have quite distinctive
properties. The field of real numbers $\br$ with its usual valuation
satisfies $|n|=n\to \infty$ for valuations of natural numbers $n$
and is said to be {\it Archimedian}. By a well known theorem of
number theory [40] the only complete Archimedian fields are those of
the real and the complex numbers. In contrast, the fields of
$p$-adic numbers, which satisfy $|n|_p \leq 1$ for all $n \in N,$
are examples of {\it non-Archimedian} fields.

Unlike the absolute value distance $|\cdot|,$ the $p$-adic
valuation satisfies the strong tringle inequality:
$$
|x+y|_p \leq \max [|x|_p,|y|_p], x, y \in \bq_p.
$$
Consequently the $p$-adic metric satisfies the strong triangle
inequality $\rho_p(x, y) \leq \max [\rho_p (x, z), \rho_p (z, y)],
x, y, z \in \bq_p,$ which means that the metric $\rho_p$ is an {\it
ultrametric,} [40]. Write $U_r(a)=\{x \in \bq_p: |x-a|_p \leq r\},$
where $r=p^n$ and $n=0, \pm 1, \pm 2, \ldots$ These are the "closed"
balls in $\bq_p$ while the sets $\bs_r(a)=\{x \in \bq_p:
|x-a|_p=r\}$ are the spheres in $\bq_p$ of such radii $r$. These
sets (balls and spheres) have a somewhat strange topological
structure from the viewpoint of our usual Euclidian intuition: they
are both open and closed at the same time, and as such are called
${\it clopen}$ sets. Finally, any $p$-adic ball $U_r(0)$ is an
additive subgroup of $\bq_p,$ while the ball $U_1(0)$ is also a
ring, which is called the {\it ring of $p$-adic integers} and is
denoted by $\bz_p.$

The $p$-adic {\it exponential function} $e^x=\sum_{n=0}^\infty
\frac{x^n}{n!}.$  The series converges in $\bq_p$ if $|x|_p \leq
r_p,$ where $r_p=1/p, p \ne 2$ and $r_2=1/4.$ $p$-adic trigonometric
functions $\sin x$ and $\cos x$ are defined by the standard power
series. These series have the same radius of convergence $r_p$ as
the exponential series.

\subsection{$p$-adic frequency model.}
As in the ordinary  probability theory [2], [3], the first $p$-adic
probability model was the frequency one, [41]. This model was based
on the simple remark that relative frequencies $\nu_N=\frac{n}{N}$
always belong to the field of rational numbers $\bq.$ And $\bq$ can
be considered as a (dense) subfield of $\br$ as well as $\bq_p$ (for
each prime number p). Therefore behaviour of sequences $\{\nu_N\}$
of (rational) relative frequencies can be studied not only with
respect to the real topology on Q, but also with respect to any
$p$-adic toplology on Q. Roughly speaking a $p$-adic probability (as
real von Mises' probability) is defined as:
\begin{equation}
\label{1} \bp (\alpha)=\lim_N \nu_N(\alpha).
\end{equation}
Here $\alpha$ is some label denoting a result of a statistical
experiment. Denote the set of all such labels by the symbol
$\Omega.$ In the simplest case $\Omega=\{0,1\}.$ Here
$\nu_{N}(\alpha)$ is the relative frequency of realization of the
label $\alpha$ in the  first $N$ trials. The $\bp (\alpha)$ is the
frequency probability of the label $\alpha.$

The main $p$-adic lesson is that it is impossible to consider, as we
did in the real case, limits of the relative frequencies $\nu_N$
when the $N \to \infty.$ Here the point "$\infty$" belongs, in fact,
to the real compactification of the set of natural numbers. So $|N|
\to \infty,$ where $|\cdot|$ is the real absolute value. The set of
natural numbers ${\bf N}$ is bounded in $\bq_p$ and it is densely
embedded into the ring of $p$-adic integers $\bz_p$ (the unit ball
of $\bq_p).$ Therefore sequences $\{N_k\}^\infty_{k=1}$ of natural
numbers can have various limits $m=\lim_{k \to \infty} N_k \in
\bz_p.$

In the $p$-adic frequency probability theory we  proceed in the
following way to provide the rigorous mathematical meaning for the
procedure (\ref{1}), see [41], [42]. We fix a $p$-adic integer $m
\in \bz_p$ and consider the class, $L_m,$ of sequences of natural
numbers $s=\{N_k\}$ such that $\lim_{k \to \infty} N_k=m$ in
$\bq_p.$

Let us consider the fixed sequence of natural numbers $s \in L_m.$
We define a $p$-adic $s$-probability as
$$
\bp(\alpha)=\lim_{k \to \infty} \nu_{N_k}(\alpha), s=\{N_k\}.
$$
This is the limit of relative frequencies with respect to the fixed
sequence $s=\{N_k\}$ of natural numbers.

For any subset $A$ of the set of labels $\Omega,$ we define its
$s$-probability as
$$
\bp (A)=\lim_{k \to \infty} \nu_{N_k}(A), s=\{ N_k \},
$$
where $\nu_{N_k}(A)$ is the relative frequency of realization of
labels $\alpha$ belonging to the set $A$ in the  first $N$ trials.

As $\bq_p$ is an additive topological semigroup (as well as $\br$),
we obtain that the $p$-adic probability is additive:

\medskip

{\bf Theorem 2.1.}
\begin{equation}
\label{2} \bp(A_1 \cup A_2)=\bp(A_1) + \bp(A_2), A_1 \cap
A_2=\emptyset.
\end{equation}

\medskip

As $\bq_p$ is even an additive topological group (as well as $\br$), we get that

{\bf Theorem 2.2.}
\begin{equation}
\label{3} \bp (A_1\setminus  A_2)=\bp(A_1)-\bp(A_1 \cap A_2).
\end{equation}

\medskip

Trivially, for any sequence $s= \{N_k\},$ $\bp(\Omega)=\lim_{k\to
\infty} \nu_{N_k}(\Omega)=1,$  as $\nu_N(\Omega)=\frac{N}{N}=1$ for
any $N.$

 As $\bq_p$ is a multiplicative
topological group (as well as $\br)$ , we obtain (see von Mises [2],
[3] for the real case and [11] for the $p$-adic case) Bayes' formula
for conditional probabilities:

\medskip

{\bf Theorem 2.3.}
\begin{equation}
\label{4} \bp (A\vert B)=\lim_{k \to \infty} \frac{\nu_{N_k}(A \cap
B)}{\nu_{N_k}(A)}=\frac{\bp(A \cap B)}{\bp(A)}, \bp(A) \ne 0.
\end{equation}

\medskip

As we know, frequency probability played the crucial role in
conventional probability theory for determination  of the range of
values (namely, the segment [0,1]) of a probabilistic measure, see
remarks on von Mises' theory in Kolmogorov's book [1]. Frequencies
always lie between zero and one. Thus their limit belongs to the
same range.

In the $p$-adic case we can proceed in the same way. Let $r \equiv
r_m = \frac{1}{|m|}_p$ (where $r=\infty$ for $m=0$). We can easily
obtain, see [42], that for the $p$-adic frequency $s$-probability,
$s \in L_m,$ the values of $\bp$ always belong to the $p$-adic ball
$U_r(0)=\{x \in \bq_p:|x|_p \leq r\}.$ In the $p$-adic probabilistic
model such a ball $U_r(0)$ plays the role of the segment [0,1] in
the real probabilistic model.

\subsection{Measure-theoretic approach}
As in the real case, the structure of an additive topological group
of $\bq_p$ induces the main properties of probability that can be
used for the axiomatization in the spirit of Kolmogorov, [1].

Let us fix $r=p^{\pm l}, l=0, 1, \ldots,$  or $r=\infty.$

{\bf Axiomatics 1.} {\it Let $\Omega$ be an arbirary set (a sample
space) and let $F$ be a field of subsets of $\Omega$ (events).
Finally, let $\bp: F \to U_r(0)$ be an additive function (measure)
such that $\bp (\Omega)=1.$ Then the triple $(\Omega, F, \bp)$ is
said to be a $p$-adic $r$-probabilistic space and $\bp\; $ $p$-adic
$r$-probability }.

 Following to Kolmogorov we should find some
technical mathematical restriction on $\bp$ that would induce
fruitful integration theory and give the possibility to define
averages. Kolmogorov (by following Borel, Lebesque, Lusin, and
Egorov) proposed to consider the $\sigma$-additivity of measures and
the $\sigma$-structure of the field of events. Unfortunately, in the
$p$-adic case the situation is not so simple as in the real one. One
could not just copy Kolmogorov's approach and consider the condition
of $\sigma$-additivity. There is, in fact, a No-Go theorem, see,
e.g., [43]:

\medskip

{\bf Theorem 2.4.} {\it All $\sigma$-additive $p$-adic valued
measures defined on $\sigma$-fields are discrete.}

\medskip

Here the difficulty is not induced by the $\sigma$-additivity, but
by an attempt to  extend a measure from the field $F$ of its
definition to a $\sigma$-field. Roughly speaking there exist
$\sigma$-additive ``continuous'' $\bq_p$-valued measures, but they
could not be extended from the field $F$ to the $\sigma$-field
generated by $F$. Therefore it is impossible to choose the
$\sigma$-additivity as the basic integration condition in the
$p$-adic probability theory.

The first important condition  (that was already invented in the
first theory of non-Archimedian integration of Monna and Springer
[44]) is {\it boundedness:}
 $$||A||_\bp=\sup\{|\bp(A)|_p: A \in F\}
< \infty.$$

Of course, if $\bp$ is a $p$-adic $r$-probability with $r<\infty,$
then this condition is fulfilled automatically. It is nontrivial
only if the range of values of  a $p$-adic probability is unbounded
in $\bq_p.$\footnote{In the frequency formalism this corresponds to
considering of $p$-adic (frequency) $s$-probabilities for $s \in
L_0;$ e.g., $s=\{N_k=p^k\}.$}. We pay attention to one important
particular case in that the condition of boundedness alone provides
the fruitful integration theory. Let $\Omega$ be a compact
zero-dimensional topological space.\footnote{There exists a basis of
neighborhoods that are open and closed at the same time.} Then the
integral
$$
E \xi=\int_\Omega \xi(\omega) \bp (d\omega)
$$
is well defined for any continuous $\xi:\Omega \to \bq_p.$ This
theory works well for the following choice: $\Omega$ is the ring of
$q$-adic integers $Z_{q},$ and $\bp$ is a bounded $p$-adic
$r$-probability, $r < \infty.$ The integral is defined as the limit
of Riemannian sums [44].

But in general boundedness alone does not imply a fruitful
integration theory. We should consider another condition, namely
{\it continuity} of $\bp,$ see appendix 3. The most general
continuity condition was proposed by A. van Rooij [43], see appendix
3.\footnote{We remark that in many cases continuity coincides with
$\sigma$-additivity, see appendix 3.}

\medskip

{\bf Definition 2.1.} {\it A $p$-adic valued measure that is
bounded, continuous, and normalized is called $p$-adic probability
measure.}

\medskip

 Everywhere below we consider
$p$-adic probability spaces endowed with $p$-adic probability
measures.

Let ($\Omega, F, \bp$) be a $p$-adic probabilistic  space. {\it
Random variables} $\xi: \Omega \to \bq_p$ are defined as
$\bp$-integrable functions, see appendix 3.

As the frequency $p$-adic probability theory induces, see [41], (as
a Theorem) Bayes' formula for conditional probability, we can use
(\ref{4}) as the definition of conditional probability in the
$p$-adic axiomatic approach (as it was done by Kolmogorov in the
real case).

{\bf Example 2.1.} ($p$-adic valued uniform distribution on the
space of $q$-adic sequences). Let $p$ and $q$ be two prime numbers.
We set $X_q=\{0, 1, \ldots, q-1\}, \Omega_q^n=\{x=(x_1, \ldots,
x_n): x_j \in X_q\}, \Omega_q^\star=\bigcup_n \Omega_q^n$ (the space
of finite sequences), and
$$
\Omega_q=\{\omega=(\omega_1, \ldots, \omega_n, \ldots): \omega_j \in X_q\}
$$
(the space of infinite sequences). For $x \in \Omega_q^n,$ we set
$l(x)=n.$ For $x \in \Omega_q^\star, l(x)=n,$ we define a cylinder
$U_x$ with the basis $x$ by $U_x=\{\omega \in \Omega_q: \omega_1 =
x_1, \ldots, \omega_n=x_n\}.$ We denote by the symbol $F_{\rm cyl}$
the field of subsets of $\Omega_q$ generated by all cylinders. In
fact, the $F_{\rm cyl}$ is the collection of all finite unions of
cylinders.

First we define the uniform distribution on cylinders  by setting
$\mu(U_x)=1/q^{l(x)}, x \in \Omega^\star_q.$ Then we extend $\mu$ by
additivity to the field $F_{\rm cyl}$. Thus $\mu: F_{\rm cyl} \to
\bq.$ The set of rational numbers can be considered as a subset of
any $\bq_p$ as well as a subset of $\br.$ Thus $\mu$ can be
considered as a $p$-adic valued measure (for any prime number $p)$
as well as the real valued measure. We use symbols $\bp_p$ and
$\bp_\infty$ to denote these measures. The probability space for the
uniform $p$-adic measure is defined as the triple
$$
{\cal P}=(\Omega, F, \bp), \;{\rm where}\; \Omega=\Omega_q, F=F_{\rm
cyl} \;{\rm and}\; \bp=\bp_p.
$$
The $\bp_p$ is called a {\it uniform $p$-adic probability
distribution.}

The uniform $p$-adic probability distribution is a  probabilistic
measure iff $p \ne q.$ The range of its values is a subset of the
unit $p$-adic ball.

{\bf Remark 2.1.} Values of $\bp_p$ on cylinders coincide with
values of the standard (real-valued) uniform probability
distribution (Bernoulli measure) $\bp_\infty.$ Let us consider, the
map $j_\infty(\omega)=\sum_{j=0}^\infty \frac{\omega_j}{2^{j+1}}.$
The $j_\infty$ maps the space $\Omega_q$ onto the segment [0, 1] of
the real line $\br$ (however, $j_\infty$ is not one to one
correspondence). The $j_\infty$-image of the Bernoulli measure is
the standard Lebesque measure on the segment [0,1] (the uniform
probability distribution on the segment [0,1]).

{\bf Remark 2.2.} The map $j_q:\Omega_q \to \bz_q,
j_q(\omega)=\sum_{j=0}^\infty \omega_j q^j,$ gives (one to one!)
correspondence between the space of all q-adic sequences $\Omega_q$
and the ring of q-adic integers $Z_q$. The field $F_{\rm cyl}$ of
cylindrical subsets of $\Omega_q$ coincides with the field
$B(\bz_q)$ of all clopen (closed and open at the same time) subsets
of $\bz_q.$ If $\Omega_q$ is realized as $\bz_q$ and $F_{\rm cyl}$
as $B(Z_q),$ then $\mu_p$ is the $p$-adic valued Haar measure on
$\bz_q$. The use of the topological structure of $\bz_q$ is very
fruitful in the integration theory (for $p \ne q).$ In fact, the
space of integrable functions $f: \bz_q \to \bq_p$ coincides with
the space of continuous functions (random variables) $C(\bz_q,
\bq_p),$ see [44], [40], [43], [11].

\section {$p$-adic limit theorems}
\subsection {$p$-adic Asymptotics of combinatorial probabilities}
Everywhere in this section $p$ is a prime number distinct from 2. We
start with considering the classical Bernoulli scheme (in the
conventional probabilistic framework) for random variables
$\xi_j(\omega)=0,1$ with probabilities $1/2, j=1,2,\ldots.$ First we
consider a finite number $n$ of random variables: $\xi_1(\omega),
\ldots, \xi_n(\omega).$ A sample space correspondding to these
random variables can be chosen as the space $\Omega_2^n =
\{0,1\}^n.$ The probability of an event $A$ is defined as
$$
\bp^{(n)}=\frac{|A|}{|\Omega_2^n|}=\frac{|A|}{2^n},
$$
where the symbol $|B|$ denotes the number of elements in a set $B$.
The typical problem of ordinary  probability theory  is to find an
asymptotic behavior of the probabilities $\bp^{(n)}(A), n \to
\infty.$ It was the starting point of the theory of limit theorems
in conventional probability theory.

But the probabilities $\bp^{(n)}(A)$ belong to the field of rational
numbers $\bq.$ We may study behavior of $\bp^{(n)}(A)$, not only
with respect to the usual real metric $\rho_\infty(x, y)$ on $\bq$,
but also with respect to an arbitrary metric $\rho(x,y)$ on $\bq.$
We have studied the case of the $p$-adic metric on $\bq,$ see [28],
[39]. We remark that $\bp^{(n)}(A)=\sum_{x \in A} \mu(U_x),$ where
$\mu$ is the uniform distribution on $\Omega_2.$ By realizing $\mu$
as the (real valued) probability distribution $\bp_\infty$ we use
the formalism of the conventional probability theory. By realizing
$\mu$ as the $p$-adic valued probability distribution $\bp_p$ we use
the formalism of $p$-adic probability theory.

What kinds of events $A$ are naturally coupled to the $p$-adic
metric? Of course, such events must depend on the prime number $p.$
As usual, we consider the sums
$$
S_n(\omega)=\sum_{k=1}^n \xi_n(\omega).
$$
We are interested in the following question. Does $p$ divide the sum
$S_n(\omega)$ or not? Set $A(p, n)=\{\omega \in \Omega_2^n: p \;
{\rm \; divides \; the \; sum } S_n(\omega)\}.$ Then $\bp^{(n)}(A(p,
n))=L(p, n)/2^n,$ where $L(p, n)$ is the number of vectors $\omega
\in \Omega_2^n$ such that $p$ divides $|\omega|=\sum_{j=1}^n
\omega_j. $ As usual, denote by $\bar A$ the complement of a set
$A.$ Thus $\bar A(p, n)$ is the set of all $\omega \in \Omega_q^n$
such that $p$ does not divide the sum $S_n(\omega).$ We shall see
that the sets $A(p, n)$ and $\bar A(p,n)$ are asymptotically
symmetric from the $p$-adic point of view:
\begin{equation}
\label{5}
\bp^{(n)}(A(p, n))\to \frac{1}{2} \;{\rm and}\; \bp^{(n)}(\bar A (p, n))\to \frac{1}{2}
\end{equation}
in the $p$-adic metric when $n \to 1$ in the same metric. Already in
this simplest case we shall see that the behavior of sums
$S_n(\omega)$ depends crucially on the choice of a sequence
$s=\{N_k\}_{k=1}^\infty$ of natural numbers. A limit distribution of
the sequence of random variables  $S_n(\omega),$ when $n \to \infty$
in the ordinary sense, does not exist. We have to describe all
limiting distributions for different sequences $s$ converging in the
$p$-adic topology.

Let ($\Omega, F, \bp$) be a $p$-adic  probabilistic space and
$\xi_n: \Omega \to \bq_p (n=1,2,\ldots)$ be a sequence of equally
distributed independent random variables, $\xi_n= 0,1$ with
probability $1/2.$\footnote{Here $1/2$ is considered as a $p$-adic
number. In the conventional theory $1/2$ is considered as a real
number.} We start with the following result that can be obtained
through purely combinatorial considerations (behavior of binomial
coefficients $C_m^r$ in the $p$-adic topology).

\medskip

{\bf Theorem 3.1.} {\it Let $m=0,1, \ldots, p^s - 1(s=1,2,\ldots), r=0, \ldots, m,$ and $l \geq s.$ Then}
$$
\lim_{n \to m} \bp (\omega: S_n(\omega) \in U_{1/p^l}(r))=\frac{C^r_m}{2^m}.
$$

\medskip

Formally this theorem can be reformulated as the following result
for the convergence of probabilistic distributions: {\it The
limiting distribution on $\bq_p$ of the sequence of the sums
$S_n(\omega)$, where $n \to m$ in $\bq_p,$ is the discrete measure}
$\kappa_{1/2, m}=2^{-m} \sum_{r=0}^m C^r_m \delta_m.$

We consider the event $A(p, n, r)=\{\omega:S_n(\omega)=pi + r\} $
for $r=0,1, \ldots, p-1.$ This event consists of all $\omega$ such
that the residue of $S_n(\omega) \mod p$ equals to $r$. Note that
the set $A(p, n, r)$ coincides with the set $\{\omega:S_n(\omega)
\in U_{1/p}(r)\}.$

\medskip

{\bf Corollary 3.1.} {\it Let $n \to m$ in $\bq_p,$ where $m=0,1,
\ldots, p-1.$ Then the probabilities $\bp^{(n)}(A(p, n, r))$
approach $C^r_m/2^m$ for all residues $r=0, \ldots, m.$}

\medskip

In particular, as $A(p, n)\equiv A(p, n, 0),$ we get (\ref {5}).
What happens in the  case $m \geq p$? We have only the following
particular result:

\medskip

{\bf Theorem 3.2.} {\it Let $n \to p$ in $\bq_p$ and $r=0, 1, 2, \ldots, p.$ Then
$$
\lim_{n \to p} \bp(\omega: S_n(\omega) \in U_{1/p^l}(r))=\frac{C^r_p}{2^p},
$$
where $s \geq 2$ for $r=0, p$ and $s\geq 1$ for $r=1, \ldots, p-1.$}

\medskip

{\bf Remark 3.1.} (Bernard-Letac asymtotics) In [45] J. Bernard and
G. Letac have studied $p$-adic asymptotic of multi binomial
coefficients. Although they did not consider the $p$-adic
probabilistic terminology (at that moment there were no physical
motivations to consider the $p$-adic generalization of probability),
their results may be interpreted as a kind of a limit theorem for
$p$-adic probability.

\subsection{ Laws of Large Numbers}
We now study
the general case of dichotomic equally distributed independent
random variables: $\xi_n(\omega)=0,1$ with probabilities $q$ and
$q^\prime=1-q, q \in \bz_p. $ We  shall study the weak convergence
of the probability distributions $\bp_{S_{N_k}}$ for the sums
$S_{N_k}(\omega).$ We consider the space $C(\bz_p, \bq_p)$ of
continuous functions $f: \bz_p \to \bq_p.$ We will be interested in
convergence of integrals
$$
\int_{\bz_p} f(x) d\bp_{S_{N_k}}(x) \to \int_{\bz_p} f(x) d\bp_S(x), f \in C(\bz_p, \bq_p),
$$
where $\bp_S$ is the limiting probability distribution (depending on
the sequence $s=\{N_k\}$). To find the limiting distribution
$\bp_S,$ we use the method of characteristic functions. We have for
characteristic functions
$$
\phi_{N_k}(z, q, a)=\int_\Omega \exp \{z S_{N_k}(\omega)\} d {\bf
P}(\omega)=(1+q^\prime(e^z-1))^{N_k}.
$$
Here $z$ belong to a sufficiently small neighborhood of zero in the
$\bq_p;$ see [25] for detail about the $p$-adic method of
characteristic functions. Let $a$ be an arbitrary number from
$\bz_p$. Let $s=\{N_k\}^\infty_{k=1}$ be a sequence of natural
numbers converging to $a$ in the $\bq_p.$ Set $\phi(z, q,
a)=(1+q^\prime(e^z-1))^a.$ This function is analytic for small $z$.
It is easy to see that the sequence of characteristic functions
$\{\phi_{N_k}(z, q, a)\}$ converges (uniformly on every ball of a
sufficiently small radius) to the function $\phi(z, q, a).$
Unfortunately, we could not prove (or disprove) a $p$-adic analogue
of Levy's theorem. Therefore in the general case the convergence of
characteristic functions does not give us anything. However, we
shall see that we have Levy's situation in the particular case under
consideration: There exists the bounded probability measure
distribution $\kappa_{q, a}$ having the characteristic function
$\phi(z, q, a)$ and, moreover, $\bp_{S_{N_k}} \to
\bp_S=\kappa_{q,a}, N_k \to a.$

We start with the first part of the above statement.  Here we shall
use Mahlers integration theory on the ring of $p$-adic integers, see
e.g., [40]. We introduce a system of binomial polynomials: $C(x,
k)=C_x^k=\frac{x(x-1)\ldots(x-k+1)}{k!}$ (that are considered as
functions from $\bz_p$ to $\bq_p$). Every function $f \in C(\bz_p,
\bq_p)$ is expanded into a series (a Mahler expansion, see [40])
$f(x)=\sum_{k=0}^\infty a_k C(x, k).$ It converges uniformly on
$\bz_p.$ If $\mu$ is a bounded measure on $\bz_p$, then
$$
\int_{\bz_p}f(x) \mu(dx)=\sum a_k \int_{\bz_p} C(x, n)\mu(dx).
$$
Therefore to define a $p$-adic valued measure on $\bz_p$ it suffices
to define coefficients $\int_{\bz_p}C(x, n)\mu(dx).$ A measure is
bounded iff these coefficients are bounded. Using the expansion of
$\phi(z, q, a),$ we obtain
$$
\la_m(q, a)=\int_{\bz_p} C(x, m)\kappa_{q, a}(dx)=(1-q)^m C(a, m).
$$
As $|C(a, m)|_p \leq 1$ for $a \in \bz_p,$ we get that the
distribution $\kappa_{q, a}$ (corresponding to $\phi(z, q, a)$) is
bounded measure on $\bz_p$. Set $\la_{mn}(q, a)=\int_\Omega
C(S_n(\omega), m) dP(\omega).$

We compute
$$\la_{mN_k}(q, a)=(1-q)^m C_{N_k}^m.$$

Thus $\la_{mN_k}(q, a) \to \la_m(q, a), N_k \to a.$ This implies the following limit theorem.

{\bf Theorem 3.3.} ($p$-adic Law of Large Numbers.) {\it The
sequence of probability  distributions $\{\bp_{S_{N_k}}\}$ converges
weakly to $\bp_S=\kappa_{q, a},$ when $N_k \to a$ in $\bq_p.$}

\subsection{The central limit theorem}
Here we restrict our considerations to the case of symmetric random
variables $\xi_n(\omega)=0, 1$ with probabilities $1/2.$ We study
the $p$-adic asymptotic of the normalized sums
\begin{equation}
\label{6}
G_n(\omega)=\frac{S_n(\omega)-ES_n(\omega)}{\sqrt{DS_n(\omega),}}
\end{equation}
Here $ES_n=n/2, D\xi_n=E\xi^2-(E\xi)^2=1/4$ and $DS_n = n/4.$ Hence
$$
G_n(\omega)=\frac{S_n(\omega)-n/2}{\sqrt{n}/2}=\sum_{j=1}^n
\frac{2\xi_n}{\sqrt {n}}-\sqrt{n.}
$$
By applying the method of  characteristic functions we can find the
characteristic function of the limiting distribution. Let us compute
the characteristic function of random variables $G_n(\omega):$
$$
\psi_n(z)=(cosh\{z/\sqrt{n}\})^n.
$$
Set $\psi(z,
a)=(cosh\{z/\sqrt{a}\})^a, a \in \bz_p, a \ne 0.$ This function
belongs to the space of locally analytic functions. There exists the
$p$-adic  analytic generalized function, see [25] for detail,
$\gamma_a$ with the Borel-Laplace transform $\psi(z,a).$
Unfortunately, we do not know so much about this distribution (an
analogue of Gaussian distribution?). We only proved the following
theorem:

{\bf Theorem 3.4.} {\it The $\gamma_1$ is the bounded measure on $\bz_p.$}

{\bf Open Problems:}

1). Boundedness of $\gamma_a$ for $a \ne 1.$

2). Weak convergence of $\bp_{G_n}$ to $\bp_G=\gamma_a$ (at least for $a=1).$

\section{Axiomatics for probability valued in a topological group}
Let $G$ be a commutative (additive) topological group. In general,
it can be nonlocally compact.\footnote{In principle, we could
proceed in the same way in the non-commutative case.}

Let us choose a fixed subset $\Delta$ of the group $G.$

{\bf Axiomatics 2.} {\it Let $\Omega$ and $F$ be as in Axiomatics 1.
Let $\bp: F \to \Delta$  be an additive function (measure). The
triple ($\Omega, F, \bp$) is said to be a $G$-probabilistic space
(with the $\Delta$-range of probability).}

We also have to add an integration condition. Such a condition
depends  on the topological structure of $G.$ It seems to be
impossible to propose a general condition providing fruitful
integration theory.

The reader might say that our definition of a $G$-probabilistic
measure is too general. Moreover,  our real probabilistic intuition
would protest against disappearence of the ${\it {unit \;
probability}}$ from consideration. We shall discuss this problem in
section 5.

We now consider a modification of the above axiomatics that includes
a kind of `unit probability'. Let $E=\bp(\Omega)$ be a nonzero
element in $G$. Let $G$ be metrizable (with the metric $\rho$). The
additional (`unit-probability') axiom should be of the following
form:
\begin{equation}
\label{7}
\sup_{A \in F} \rho(0, \bp(A))=\rho(0, E).
\end{equation}
A $G$-probabilistic space in that (\ref{7}) holds true is called a
$G$-probabilistic space with {\it unit probability axiom.} Of
course, the consideration of such probabilistic spaces seems to be
more natural from the standard probabilistic viewpoint. Therefore it
would be natural to start with consideration of such models.
However, for many important $G$-probabilistic spaces the unit
probability axiom does not hold true. At the moment we know a few
examples of $G$-probabilities having applications:

1). $G=\br$ and $\Delta=[0,1]$(the conventional probability theory);

2). $G=\br$ and $\Delta=\br$ (`negative probabilities', see e.g.,
[14], [18], [11], they are realized as signed measures, charges).

3). $G= {\bf C}$ and $\Delta= {\bf C}$ (`complex probabilities', see
e.g., [19], [20], they are realized as ${\bf C}$-valued measures).

4). $G=\bq_p$ and $\Delta$ is a ball in $\bq_p$ (`$p$-adic
probabilities', [42], [11], they are realized as $\bq_p$-valued
measures).

The $p$-adic model can be essentially (and rather easily)
generalized. Let  $K$ be an arbitrary complete non-Archimedian field
with the valuation (absolute value) $|\cdot|.$ We can define
$K$-valued probabilistic measures by using the same integration
conditions as in the $p$-adic case, namely boundedness and
continuity, see appendix 3.

We note that in all considered examples the  additive group $G$ has
the additional algebraic structure, namely the field structure.  The
presence of such a field structure gives the possibility to develop
essentially richer probabilistic calculus than in the general case.
Here we can introduce conditional probability by using Bayes'
formula and define the notion of independence of events.

The following slight generalization gives the  possibility to
consider a few new examples. Let $G$ be a non-Archimedian normed
ring. To simplify considerations, we again consider the commutative
case. Here:

(1) $||x|| \geq 0, ||x||=0 \leftrightarrow x= 0; $

(2)  $||x||\;||y|| \leq ||x||\;||y||$ and $||x + y|| \leq \max
(||x||, ||y||).$

We set, for $A \in F, ||A||_\bp=\sup\{||\bp(B)||:B \in F, B \subset
A\}.$ We define a $G$-probabilistic measure as a normalized
$G$-valued measure satisfying to the conditions of boundedness and
continuity, compare with appendix 3. Corresponding integration
theory is developed in the same way as in the  case of a
non-Archimedian field. One of the most important examples of
non-Archimedian normed rings is a ring of $m$-adic numbers $\bq_m,$
where $m \ne p^k, p-prime.$ It is a locally compact ring. We can
present numerous examples of non-Archimedian normed rings by
considering various functional spaces of $\bq_p$ (or $\bq_m$)-valued
functions.

For a ring $G$, we can define averages  for $G$-valued random
variables, $\xi: \Omega \to G.$ In particular, we can represent the
probability distribution of the sum $\eta=\xi_1 + \xi_2$ of two
$G$-valued random variables as the convolution of corresponding
probability distributions. Here we define the convolution of two
$G$-valued measures on $G$ as:
$$
\int_G f(x) M_1 \star M_2(dx)=\int_{G \times G} f(x_1 + x_2) M_1(dx_1) M_2(dx_2),
$$
where $f:G \to G$ is a ``sufficiently good'' function. If $G$ is a
ring and $A \in F$ is such that $\bp(A)$ is invertible, then we can
define conditional probabilities by using Bayes' formula.

We obtain a large class of new mathematical problems related to
$G$-probabilistic models. We emphasize that, despite a rather common
opinion, the probability theory is not just a part of functional
analysis (measure theory). Probability theory has also its own
ideology. The probabilistic ideology induces its own problems. Such
problems would be impossible to formulate in the framework of
functional analysis (of course, methods of functional analysis can
be essentially used for the investigation of these problems).

One of the most important problems is to find analogues of limit
theorems, compare, e.g., with [46].

{\bf Open problem:}

Let $s=\{N_k\}$ be a sequence of natural numbers and let
$$M_{11}, M_{12}, M_{1n_1};$$
$$M_{21}, M_{22}, \ldots, M_{2n_2};$$
$$M_{k1}, M_{k2}, \ldots, M_{kN_k};$$
be $G$-probabilistic measures. As usual, we have to study behavior
of convolutions:

$\alpha_k=M_{k1} \star M_{k2} \star \ldots \star M_{kN_k}$

to find analogues of limit theorems. For example, an analogue of the
law of large numbers could be formulated in the  following way. Let
$d$ be a nonzero element of a topological additive group $G$. Let
$s=\{N_k\}^\infty_{k=1}$ be a sequence of natural numbers. Suppose
that the corresponding sequence $\{N_k d\}^\infty_{k=1}$ of elements
of $G$ converges to some element $a \in G$ or to $a=\infty.$ The
latter has the standard meaning: for each neighborhood $U$ of zero
in $G$ there exists $N$ such that $N_kd \not\in U$ for all $N_k \geq
N.$

Let $(\Omega, F, \bp)$ be a $G$-probabilistic space.  Let
$\xi_n(\omega)=0, d,$ with $G$-probabilities $q, q^\prime=E-q$ where
$E=\bp(\Omega),$ be a sequence of independent random variables. Let
$S_n(\omega)$ be the sum of $n$ first variables.

{\bf Open Problem:}

{\it Does the sequence of probability distributions $\bp_{S_{N_k}}$
converge weakly to some  probability distribution $\bp_S,$ when $k
\to \infty$?}

The simplest variant of this problem is to generalize Theorem 3.1:
to find (if it exists) $\lim_{k \to \infty} \bp (S_{N_k}(\omega) \in
U_r (0)).$ In the case of a metrizable group $G,$  $U_r(0)=\{g \in
G: \rho(g,0) \leq r\}, r > 0,$ is a ball in $G$. In the case when
$G$ is a field we can consider normalized sums (\ref{6}) and try to
get an analogue of the central limit theorem.

{\bf Remark 4.1.} The presented axiomatics of the group valued
probabilities  was strongly based on authors investigations on the
$p$-adic probability theory. By using the $p$-adic experience we
understood how $G$-valued probabilities `must look'. The presented
interpretation of $G$-valued probabilities was developed in the
following way. By studying $\bq_p$-valued probabilities we
understood that we could not more use the order relation between
probabilities. In  conventional probability theory we can say, e.g.,
that $\bp(A)=1/2$ is less than $\bp(B)=2/3.$ Therefore the
occurrence of the event $A$ is less favorable than the occurrence of
the event $B$. On the other hand, there is no order structure on
$\bq_p.$ We could not say that $A$ with $\bp(A)=1/2 \in \bq_p$ is
less (or more) favorable than $B$ with $\bp(B)=2/3 \in \bq_p.$ There
is no order structure on the set of $p$-adic
probabilities\footnote{We would like to remark that the role of the
order structure on the set [0,1] of conventional probabilities is
overestimated. For example, let $\bp(A)=11/17$ and $\bp(B)=13/19.$
Of course, $\bp(A)<\bp(B)$ in $\br$. However: are you sure that you
will be rich if you play in the favor of $B$?}. It is possible to
introduce a partial order structure [51] on $\bq_p.$ However, the
formal use of that partial order structure for $p$-adic
probabilities has unexpected implications, [51]. Practically
impossible as well as practically definite events are not more
coupled to two isolated points in the set of probabilities, namely
to $\bp=0$ and $\bp=1$ respectively. Practically impossible and
definite probabilities are represented by some sets $U_0$ (such that
$0 \in U_0)$ and $U_1$ (such that $1 \in U_1),$ respectively, of
$p$-adic numbers, [51]. In the $p$-adic framework a set $A$ is
practically impossible if $\bp(A) \in U_0$ and practically definite
if $\bp(A) \in U_1.$

This $p$-adic construction  was one of motivations (at least
mathematical) of the theory of abstract models of probability
proposed by V. Maximov [52]. Maximov's formalism is based on the
deep investigation on the role of certainty, uncertainty and
randomness in probability theory. His main idea was that to apply
probability one should introduce  only the notions of certainty and
uncertainty. They can be described by some sets $U_0$ and $U_1$.
Therefore the order structure on the set of probabilities does not
play any role. In this way it is possible to proceed to
probabilities valued in abstract sets (with a rather complex
semi-algebraic structure that should reproduce the main relations
between real probabilities, see [52]). The abstract probabilistic
model of V. Maximov played the stimulating role in my investigations
to generalize the $p$-adic probability theory to an arbitrary
topological group. However, I do not support Maximov's viewpoint to
the role of certainty and uncertainty in applications of probability
theory. Each experiment must induce its special level of statistical
significance. In ordinary statistics this level is given by the
$\epsilon$-levels of statistical significance, in fact, by the
$\epsilon$-neighborhood of $\bp=0.$ Our idea was to consider an
additive topological group $G$ and probabilities valued in some
subset $U \in G, 0 \in U,$ and levels of significance are given by
neighborhoods $V$ (in the metric case $V_\epsilon, \epsilon > 0$) of
zero, see the next section for the details.

\section{Statistical interpretation of probabilities with values in a topological group}
In fact, Kolmogorov's probability theory has two (more or less
independent) counterparts: (a) axiomatics (a mathematical
representation); (i) interpretation (rules for application). The
first part is the measure-theoretic formalism. The second part is a
mixture of frequency and ensemble interpretations: {\small {"... we
may assume that to an event $A$ which has the following
characteristics: (a) one can be practically certain that if the
complex of conditions $\sum$ is repeated a large number of times,
$N$, then if $n$ be the number of occurrences of event $A$, the
ratio $n/N$ will differ very slightly from $\bp(A);$ (b) if $\bp(A)$
is very small, one can be practically certain that when conditions
$\sum$ are realized only once the event $A$ would not occur at all",
[1].}}

As we have already noticed, (a) and (i) are more or  less
independent. Therefore the Kolmogorov's measure-theoretic formalism,
(a), is used successfully, for example, in subjective probability
theory.

In practice we apply Kolmogorov's (conventional) interpretation,
(i), in the following way. First we have to fix $0 < \epsilon < 1,$
{\it significance level}. If the probability $\bp(A)$ of some events
$A$ is less than $\epsilon,$ this event is considered as practically
impossible.

It is already evident  how to generalize the conventional
interpretation of probability to $G$-valued probabilities. First we
have to fix some neighborhood of zero, $V$, a {\it significance
neighborhood}.

If the probability $\bp(A)$ of some  event $A$ belongs to $V$, this
event is considered as practically impossible.

If a group $G$ is metrizable, then  the situation is even more
similar to the standard (real) probability. We choose $\epsilon > 0$
and consider the ball $V_\epsilon=\{x \in G: \rho(0, x)<
\epsilon\}.$ If $\rho(0, \bp(A)) < \epsilon,$ then the event $A$ is
considered as practically impossible.

Let us borrow some ideas from statistics.  We are given a certain
sample space $\Omega$ with an associated distribution $\bp.$ Given
an element $\omega \in \Omega,$ we want to test the hypothesis
``$\omega$ belongs to some reasonable majority.'' A reasonable
majority $\cm$ can be described by presenting {\it critical regions
$\Omega^{(\epsilon)}(\in F)$} of the significance level $\epsilon, 0
< \epsilon < 1:\bp (\Omega^{(\epsilon)})<\epsilon.$ The complement
$\bar \Omega^{(\epsilon)}$ of a critical region
$\Omega^{(\epsilon)}$ is called $(1-\epsilon)$ confidence interval.
If $\omega \in \Omega^{(\epsilon)},$ then the hypothesis ``$\omega$
belongs to majority $\cm$'' is rejected with the significance level
$\epsilon.$ We can say that $\omega$ fails the test to belong to
$\cm$ at the level of critical region $\Omega^{(\epsilon)}.$

$G$-statistical machinery  works in the same way. The only
difference is that, instead of significance levels $\epsilon$, given
by real numbers, we consider significance levels $V$ given by
neighborhoods of zero in $G.$ Thus we consider critical regions
$\Omega^{(V)}(\in F):$
$$
\bp(\Omega^{(V)}) \in V.
$$
If $\omega \in \Omega^{(V)},$ then the hypothesis ``$\omega$ belongs
to majority $\cm$'' (represented by the statistical test
$\{\Omega^{(V)}\}$) is rejected with the significance level $V$. If
$G$ is metrizable, then we have even more similarity with the
standard (real) statistics. Here $V=V_\epsilon, \epsilon > 0.$

Of  course, the strict mathematical description of the above
statistical considerations can be presented in the framework of
Martin-Löf [9], [10], [7] statistical tests. We remark that such a
$p$-adic framework was already developed in [11]. In the $p$-adic
case (as in the real case) it is possible to enumerate effectively
all $p$-adic tests for randomness. However, a universal $p$-adic
test for randomness does not exist [11]. If the group $G$ is
metrizable we can proceed in the same way as in the real and
$p$-adic case [11] and define $G$-random sequences, namely sequences
$\omega=(\omega_1, \ldots, \omega_N, \ldots), \omega_j=0,1,$ that
are random with respect to a $G$-valued probability distribution.
However, if $G$ is not metrizable, then the notion of a recursively
enumerable set would not be more the appropriative basis for such a
theory. In any case we have an interesting

{\bf Open problem:}

{\it Development of  randomness  theory for an arbitrary topological
group.}

The general scheme of the application of $G$-valued probabilities is
the same as in the ordinary case:

1) we find initial probabilities;

2) then we perform calculations by using calculus of $G$-valued
probabilities;

3) finally, we apply the above interpretation to resulting
probabilities.

\medskip

Of course,  the main question is ``How can we find initial
probabilities?'' The situation here is more or less similar to the
situation in the ordinary probability theory. One of possibilities
is to apply the frequency arguments (as R. von Mises). We have
already discussed such an approach for $p$-adic probabilities, see
also further considerations in section 6. Another possibility is to
use subjective approach to probability. I think that everybody
agrees that there is nothing special in segment $[0,1]$ as the set
of labels for the measure of belief in the occurrence of some event.
In the same way we can use, for example, the segment [-1,1] (signed
probability) or the unit complex disk (complex probability) or the
set of $p$-adic integers $\bz_p$ ($p$-adic probability). If $G$ is a
field we can apply the machinery of Bayesian probabilities and,
finally, use our interpretation of probabilities to make a
probabilistic (`statistical') decision. The third possibility is to
use symmetry arguments, `Laplacian approach'. For example, by such
arguments we can choose (in some situations) the uniform
$\bq_p$-valued distribution.

We now turn back to the role of the unit probability and,  in
particular, axiom (\ref{7}). In fact, by considering the
interpretation of probability  based on the notion of the
significance level we need not pay the special attention to the
probability $E=\bp(\Omega).$ It is enough to consider $V$-impossible
events, $V=V(0).$ If $V$ is quite large and $\bp(A)\not \in V,$ then
an event $A$ can be considered as practically definite.

{\bf Example 5.1.} (A $p$-adic statistical test) Theorem 3.1.
implies that, for each $p$-adic sphere ${\bf S}_{1/p^l}(r),$ where
$l, r, m$ were done in Theorem 3.1:
$$
\lim_{k \to \infty} \bp (\{\om \in \Omega_2: S_{N_k}(\omega) \in
{\bf S}_{1/p^l}(r) \})=0,
$$
for each sequence $s=\{N_k\}, N_k \to m, k \to \infty.$ We can
construct a statistical test on the basis of this limit theorem (as
well as any other limit theorem). Let $s=\{N_k\}, N_k \to m,$ be a
fixed sequence of natural numbers. For any $\epsilon > 0,$ there
exists $k_\epsilon$ such that, for all $k \geq k_\epsilon,$
$$
|\bp (\{\omega \in \Omega_2: S_{N_k}(\omega) \in {\bf
S}_{1/p^l}(r)\})|_p < \epsilon.
$$
We set $\Omega^{(\epsilon)}=\bigcup_{k \geq k_\epsilon}\{\omega \in
\Omega_2: S_{N_k}(\omega) \in {\bf S}_{1/p^l}(r)\}.$ We remark that
$$
|\bp(\Omega^{(\epsilon)})|_p < \epsilon.
$$
We now define reasonable majority of outcomes as sequences that do
not belong to the sphere ${\bf S}_{\frac{1}{p^l}}(r),$
``nonspherical majority.'' Here the set $\Omega^{(\epsilon)}$ is the
critical region on the significance level $\epsilon.$

Suppose that a sequence $\omega$ belongs  to the set
$\Omega^{(\epsilon)}.$ Then the hypothesis ``$\omega$ belongs
nonspherical majority'' must be rejected with the significance level
$\epsilon.$ In particular, such a sequence $\omega$ is not random
with respect to the uniform $p$-adic distribution on $\Omega_2.$ If
, for some sequence of 0 and 1, $\omega=(\omega_j)$ we have
$\omega_1 + \ldots + \omega_{N_k} - r=\alpha$ $\; \mod p^l,
\alpha=1, \ldots, p-1,$ for all $k \geq k_\epsilon,$ then it is
rejected.

The simplest test is given by $m=1, r=0, N_k=1+p^k$  and $\omega_1 +
\ldots + \omega_{N_k}= \alpha$ $\; \mod p, \alpha=1, \ldots, p-1.$

\section{Generalized frequency models}
We proposed a theory of $G$-valued probabilities in the
measure-theoretic framework.\footnote{However, we use the frequency
experience of the real and $p$-adic probabilities to find reasonable
properties of `probability'.} In principle, we can also proceed in
the frequency framework. However, an arbitrary topological group is
too general base for such frequency probabilities. We have to start
with a topological field $T$ that contains the field of rational
numbers as a dense subfield. We proceed in the same way as in the
$p$-adic case. Let $s=\{N_k\}$ be a sequence of natural numbers
converging in $T$ to some $m$. A $T$-valued $s$-probability is
defined as the limit  $\bp=\lim_{k \to V_{N_k}} \in T$ (if it
exists).

As $T$ is  a topological field, we get additivity (\ref{2}), formula
(\ref{3}) and Bayes' formula (\ref{4}). The range of values of such
a frequency probability depends on the sequence $s$ and the topology
on $T$.

In fact, Kolmogorov's theory does not look extremely
anti-frequencist due to the presence of the strong law of large
numbers. Of course, this law was strongly criticized from the
frequency point of view, see e.g., von Mises [2], [3]. The main
critical argument is that we could not say anything about behaviour
of frequencies for a concrete sequence of trials. Nevertheless,
there are no problems in average. Therefore, by obtaining a kind of
law of large numbers (of course, in the case when the field of
rational numbers $Q$ is dense in $T$) we could strongly improve the
measure-theoretic approach for $T$-valued probabilities. A kind of
such a law we have in the $p$-adic case.

\section{Appendixes}
\subsection{Conventional stochastic processes in
$p$-adic physics}

 We consider the A-model, see section 1. Thus the
variable $x$ is $p$-adic, but the wave function $\psi$ is complex
valued. Here, as $|\psi(x)|^2 \in \br_+,$ we can apply Kolmogorov's
probability theory. However, the special structure of the sample
space $\Omega=\bq_p$ - the field of $p$-adic numbers, gives the
possibility to provide more deep probabilistic investigations. Of
course, $\bq_p$ is a locally compact additive group. Here we can
apply the general theory of probability on locally compact groups,
see, e.g., [46]-[50]. However, the importance of this concrete
example of a locally compact group was not well recognized in the
community - `Probability on Topological Structures'. On the other
hand, from the physical point of view the group $\bq_p$ is
practically not less important than the group $\br.$

Moreover,  members of the community - `Probability on Topological
Structures' were bounded by knowing the following $p$-adic No-Go
theorem:

{\bf Theorem 7.1.} {\it There are no Gaussian measures on $\bq_p$.}

This fact can easily be extracted from Heyer's book [46].  The
$\bq_p$ is a totally disconnected locally compact topological group.
This result can satisfy every mathematician. However, physicist
would not be happy! There exists a differential operator of the
second order, Vladimirov-Laplace operator, see [24], that plays the
fundamental role in $p$-adic string theory and quantum mechanics.
Physicists describe important dynamics by using  this operator. It
seems that there should be probabilistic description of these
evolutions. Of course, the corresponding stochastic processes (if it
exists at all) could not be a process with continuous trajectories.
We have to modify Bendikov's program, see e.g., [53] of the
construction of Brownian motions on locally compact groups. It seems
that a generalization of Bendikov's program for totally disconnected
locally compact groups must be based on the description of
generators of stochastic processes as second order Vladimirov-type
differential operators. It may occur that the corresponding
stochastic process will be a generalized stochastic process.

On the other hand, there were performed extended investigations on
limit theorems for general locally compact groups, see, e.g., [46].
The most interesting for us are investigations for disconnected
groups (in particular, $p$-adic), see e.g., [47]-[50]. These
investigations can be of great importance for $p$-adic physics and
cognitive sciences.

\subsection{Limit theorem for signed, complex and Banach-algebras
valued probabilities}

It seems that  the first limit theorem for signed probabilities was
considered in the paper of [54]. Then such probabilistic
distributions were considered by V. Krylov [55] in the connection to
the Cauchy problem for the equations:
\begin{equation}
\label{8}
\frac{\delta f}{\delta t}=(-1)^{q+1} \frac{\delta^{2q}f}{\delta x^{2q}}
\end{equation}
He proved a limit theorem for the fundamental solution of (\ref{8})
that generalized the standard central limit theorem (for $q=1$).

The central limit theorem for Schr\"odinger's equation was proved by
A. Khrennikov and O. Smolyanov [56]. Then it was generalized to
arbitrary complex (in particular, signed) distributions (on finite
as well as infinite dimensional locally convex spaces), see [57]. We
also have to mention an extended research activity on limit theorems
for fundamental and pseudo-differential equations. However, the
author is not a specialist in this field. Therefore we will not
concern these investigations in the present review. Finally, we
remark that `probability measures' with values in supercommutative
Banach superalgebras (in particular, commutative) were used in the
framework of so called supersymmetric physics [58]. Here I proved an
analogue of the central limit theorem for Gaussian and Feynman
superdistributions.

\subsection{Measures that take values in non-Archimedian fields}

Let $X$ be an arbitrary set and let ${\cal R}$ be a ring  of subsets
of $X$. The pair (${\cal X, R}$) is called a {\it measurable space.}
The ring ${\cal R}$ is said to be {\it separating} if for every two
distinct elements, $x$ and $y$, of $X$ there exists an $A \in {\cal
R}$ such that $x \in A, y \not \in A.$ We shall consider measurable
spaces only over separating rings which cover $X$.

A subcollection ${\cal S}$ of ${\cal R}$ is said to be {\it
shrinking}  if the intersection of any two elements of ${\cal S}$
contains an element of ${\cal S}$. If ${\cal S}$ is shrinking, and
if $f$ is a map ${\cal R} \to K$ or ${\cal R} \to \br,$ we say that
$\lim_{A \in {\cal S}} f(A)=0$ if for every $\epsilon > 0,$ there
exists an $A_0 \in {\cal S}$ such that $|f(A)|\leq \epsilon$ for all
$A \in {\cal S}, A \subset A_0.$

Let $K$ be a  non-Archimedian field with the valuation $|\cdot|.$ A
{\it measure} on ${\cal R}$ is a map $\mu {\cal R} \to K$ with the
properties:

(i) $\mu$ is additive;  (ii) for all $A \in {\cal R},
||A||_\mu=\sup\{|\mu(B)|:B \in {\cal R}, B \subset A\} < \infty;
(iii)$ if $\cs \subset {\cal R}$ is shrinking and has empty
intersection, then $\lim_{A \in \cs}\mu(A)=0.$

We call these conditions respectively {\it additivity, boundedness, continuity.} The latter condition is equivalent to the following: $\lim_{a \in \cs}||A||_\mu=0$ for every shrinking collection $\cs$ with empty intersection.

Condition (iii) is the replacement for $\sigma$-additivity. Clearly (iii) implies $\sigma$-additivity. Moreover, for the most interesting cases, (iii) is equivalent to $\sigma$-additivity [43]. Of course, we could in principle restrict our attention to these cases and use the standard condition of $\sigma$-additivity. However, in that case we should use some topological restriction on the space $X$. This implies that we must consider some topological structure on a $p$-adic probability space. We do not like to do this. We would like to develop the theory of $p$-adic probability measures in the same way as A. N. Kolmogorov (1933) developed the theory of real valued probability measures by starting with an arbitrary set algebra.

For any set $D$, we denote its characteristic function by the symbol $I_D.$ For $f:X \to K$ and $\phi:X \to [0, \infty),$ put

$$||f||_\phi=\sup_{x \in X}|f(x)|\phi(x).$$

We set

$$N_\mu(x)=\inf_{U \in {\cal R}, x \in U}||U||_\mu$$
for $x \in X.$ Then $||A||_\mu=||I_A||_{N_\mu}$ for any $A \in {\cal
R}.$ We set $||f||_\mu=||f||_{N_\mu}.$ A {\it step function} (or
${\cal R}$-step function) is a function $f:X \to K$ of the form
$f(x)=\sum_{k=1}^N c_k I_{A_k}(x)$ where $c_k \in K$ and $A_k \in
{\cal R}, A_k \cap A_1=\emptyset, k \ne l.$ We set for such a
function

$$\int_X f(x) \mu (dx)=\sum_{k=1}^N c_k \mu(A_k).$$

Denote the space of all step functions by the symbol $S(X).$  The
integral $f \to \int_X f(x) \mu (dx)$ is the linear functional on
$S(X)$ which satisfies the inequality $|\int_X f(x) \mu
(dx)|\leq||f||_\mu.$ A function $f:X \to K$ is called $\mu$-{\it
integrable} if there exists a sequence of step functions $\{f_n\}$
such that $\lim_{n \to \infty}||f-f_n||_\mu=0.$ The $\mu$-integrable
functions form a vector space $L(X, \mu)$ (and $S(X) \subset L(X,
\mu)$). The integral is extended from $S(X)$ on $L(X, \mu)$ by
continuity. Let ${\cal R}_\mu=\{A: A \subset  X, I_A \in L(X,
\mu)\}.$ This is a ring. Elements of this ring are called
$\mu$-measurable sets. By setting $\mu(A)=\int_X I_A(x) \mu (dx)$
the measure $\mu$ is extended to a measure on ${\cal R}_\mu.$ This
is the {\it maximal extension} of $\mu,$ i.e., if we repeat the
previous procedure starting with the ring ${\cal R}_\mu,$ we will
obtain this ring again. Finally, we mention  investigation for
solenoids, see e.g., [53] applications of $p$-adic valued measures
in number theory, see e.g., [40] and recent paper of D.
Neuenschwander [59] on $p$-adic valued measures.

\medskip

{\bf Acknowledgements:}

I would like to thank L. Accardi, S. Albeverio, A. Bendikov,  A.
Bulinskii, W. Hazod, H. Heyer, T. Hida, S. Kozyrev, G. Letac, V.
Maximov, D. Neuenschwander, A. Shiryaev, Yu. Prohorov,  O.
Smolyanov, V. Vladimirov, and I. Volovich for fruitful discussions.
I am especially thankful to A. Bulinskii and A. Shiryaev for
conversations on original view of A. N. Kolmogorov on axiomatization
of probability and to A. Accardi and V. Maximov for many years of
debates on non-Kolmogorovian extensions of probability theory. I
also would like to thank V. P. Maslov for a series of discussions on
frequency probability and quantum probability as well as their
applications outside quantum domain and in particular in financial
mathematics and his support of my frequency analysis of foundations
of probability theory.

\medskip

{\bf References}

[1] A. N. Kolmogorov, {\it Grundbegriffe der Wahrscheinlichkeitsrechnung.}
Springer Verlag, Berlin (1933); reprinted: {\it Foundations of the Probability Theory.} Chelsey Publ. Comp., New York (1956).

[2] R. von Mises, Grundlagen der Wahrschienlichkeitsrechnung. {\it Math. Z.,} {\bf 5}, 52-99 (1919).

[3] R. von Mises, {\it The mathematical theory of probability and statistics} (Academic , London, 1964)

[4] E. Kamke, Über neuere Begründungen der Wahrscheinlichkeitsrechning. {\it Jahresbericht der Deutschen Matematiker-Vereinigung,} {\bf 42}, 14-27 (1932).

[5] J. Ville, {\it Etude critique de la notion de collective.} Gauthier- Villars, Paris (1939).

[6] E. Tornier, {\it Wahrscheinlichkeitsrechnung und allgemeine Integrationsteorie.,} Univ. Press, Leipzing (1936).

[7] M. Li M., P. Vitànyi, {\it An introduction to Kolmogorov
complexity and its applications.}  Springer, Berlin-Heidelberg-New
York (1997).

M. Van Lambalgen, Von Mises' definition of random sequences
reconsidered. {\it J. of Symbolic Logic,} {\bf 52,} N. 3 (1987).

[8] A. N. Kolmogorov, Three approaches to the quantitative
definition of information.  {\it Problems Inform. Transmition.,}
{\bf 1}, 1-7 (1965).

 A. K. Zvonkin and L. A. Levin, A complexity of finite objects and
justification of the notions of information and randomness on the
basis of algorithms.{\it Uspehi Matem. Nauk}, {\bf 25}, 85--127
(1970).

[9] P. Martin-Löf, On the concept of random sequences. {\it Theory of Probability Appl.,} {\bf 11}, 177-179 (1966).

[10] P. Martin-Löf, The definition of random sequences. {\it Inform. Contr.,} {\bf 9}, 602-619 (1966).

[11] A. Yu. Khrennikov, {\it Interpretation of probability.} (VSP Int. Publ., Utrecht. 1999).

[12] A. Yu. Khrennikov, A perturbation of CHSH inequality induced by fluctuations of ensemble distributions. {\it J. of Math. Physics,} {\bf 41}, N.9, 5934-5944 (2000).

[13] A. Yu. Khrennikov, Statistical measures of ensemble nonreproducibility and correction to Bell's inequality. {\it Il Nuovo Cimento,} {\bf B 115}, N.2, 179-184 (2000).

[14] P. A. M. Dirac, The physical interpretation of quantum mechanics. {\it Proc. Roy. Soc. London,} {\bf A 180}, 1-39 (1942).

[15] E. Wigner, {\it Quantum-mechanical distribution functions revisted,} in: {\it Perspectives in quantum theory}. Yourgrau W. and van der Merwe A., editors, MIT Press, Cambridge MA (1971).

[16] R. P. Feynman, Negative probability. {\it Quantum Implications, Essays in Honour of David Bohm,} B. J. Hiley and F. D. Peat, editors. Routledge and Kegan Paul, London, 235-246 (1987).

[17] W. Muckenheim, A review on extended probabilities. {\it Phys. Reports,} {\bf 133}, 338-401 (1986).

[18] A. Yu. Khrennikov, $p$-adic description of Dirac's hypothetical world with negative probabilities. {\it Int. J. Theor. Phys.,} {\bf 34}, No. 12, 2423-2434 (1995).

[19] P. A. M. Dirac, On the analogy between classical and quantum mechanics. {\it Rev. of Modern Phys.,} {\bf 17}, N 2/3, 195-199 (1945).

[20] E. Prugovecki, Simultaneous measurement of several observables. {\it Found. of Physics,} {\bf 3}, No. 1, 3-18 (1973).

[21] Volovich I. V., $p$-adic string. {\it Class. Quant. Grav.,} {\bf 4}, 83-87 (1987).

[22] P. G. O. Freund and M. Olson, Non-Archimedian strings. {\it Phys. Lett. B.} {\bf 199}, 186-190 (1987).

[23] P. G. O. Freund and E. Witten, Adelic string amplitudes. {\it Phys. Lett. B.,} {\bf 199}, 191-195 (1987).

[24] V. S. Vladimirov, I. Volovich, and E. I. Zelenov, {\it $p$-adic analysis and Mathematical Physics.} World Scientific Publ. Singapore (1993).

[25] A. Yu. Khrennikov, {\it $p$-adic valued distributions and their applications to the mathematical physics,} Kluwer Acad. Publishers, Dordrecht (1994).

[26] A. Yu. Khrennikov, {\it Non-Archimedian analysis: quantum
paradoxes, dynamical systems and biological models.} Kluwer Acad.
Publishers, Dordrecht (1997).

[27] S. Albeverio, A. Yu. Khrennikov, Representation of the Weyl group in spaces of square integrable functions with respect to $p$-adic valued Gaussian distributions. {\it J. of Phys. A, 29}, 5515-5527 (1996).

[28] A. Yu. Khrennikov, $p$-adic probability distribution of hidden variables, {\it Physica A}, {\bf 215}, 577-587 (1995).

[29] A. Yu. Khrennikov, $p$-adic probability interpretation of Bell's inequality paradoxes. {\it Physics Letters A,} {\bf 200}, 119-223 (1995).

[30] E. Thiran E., Verstegen D., Weyers J., $p$-adic dynamics, {\it J. Stat. Phys.} {\bf 54} (1989), 893-913.

[31] J. Lubin, Non-archimedian dynamical systems, {\it Compositio Mathematica,} {\bf 94}, 1994, 321-346.

[32] A. Yu. Khrennikov, M. Nilsson, {\it $p$-adic deterministic and
random dynamical systems.} Kluwer, Dordreht (2004).

S. Albeverio, M. Gundlach,  A. Yu. Khrennikov, K.- O. Lindahl, On
Markovian behaviour of $p$-adic random dynamical systems. {\it
Russian J. of Math. Phys.}, {\bf 8} (2), 135-152 (2001).

A. Yu. Khrennikov, M. Nilsson, On the number of cycles for $p$-adic
dynamical systems. {\it J. Number Theory,} {\bf 90}, 255-264 (2001).

[33] S. Albeverio, A. Yu. Khrennikov,  P. Kloeden, Memory retrieval
as a $p$-adic dynamical system. {\it Biosystems,} 49, 105-115
(1999).

A. Yu. Khrennikov, $p$-adic discrete dynamical systems and their
applications in physics and cognitive sciences. {\it Russian J.
Math. Phys.}, {\bf 11,} N 1, 45--70 (2004).

8. Khrennikov A.Yu.,{\it Information dynamics in cognitive,
psychological and anomalous phenomena.} Kluwer, Dordreht (2004).

[34] T. L. Fine, {\it Theories of probabilities, an examination of foundations.} Academic Press, New York (1973).

[35] Accardi L., The probabilistic roots of the quantum mechanical paradoxes. {\it The wave-particle dualism. A tribute to Louis de Broglie on his 90th Birthday,} Edited by S. Diner, D. Fargue, G. Lochak and F. Selleri, D. Reidel. Publ. Company, Dordrecht, 47-55 (1970).

[36] S. Gudder, Probability manifolds. {\it J. Math. Phys.,} {\bf 25}, 2397-2401 (1984).

[37] {\it Selected works of A. N. Kolmogorov,} {\bf 2}, {\it Probability theory and mathematical statistics.} Ed. A. N. Shiryayev. Kluwer Academic Publ., Dordrecht (1986).

[38] A. Yu. Khrennikov, Laws of large numbers in non-Archimedian probability theory. {\it Izvestia Akademii Nauk,} {\bf 64}, N.1, 211-223 (2000).

A. Yu. Khrennikov, Limit behaviour of sums of independent random
variables with respect to the uniform $p$-adic distribution. {\it
Statistics and Probability Lett.,} {\bf 51}, 269-276 (2001).

[39]  V. P. Maslov, {\it Quantization and ultra-second
quantization.} Institute of Computer Investigations, Moscow (2001).

[40] W. Schikhov, {\it Ultrametric calculus.} Cambridge Univ. Press, Cambridge (1984).

[41] A. Yu. Khrennikov, $p$-adic probability and statistics. {\it Dokl. Akad. Nauk USSR,} {\bf 322,} 1075-1079 (1992).

[42] A. Yu. Khrennikov, An extension of the frequency approach of R. von Mises and the axiomatic approach of N. A. Kolmogorov to the $p$-adic theory of probability. {\it Theory of Probability and Appl., }{\bf 40}, N. 2, p. 458-463 (1995).

[43] A. Van Rooij, {\it Non-archimedian functional analysis.} Marcel Dekker, Inc., New York (1978).

[44] A. Monna and T. Springer, Integration non-Archimédienne, 1, 2. {\it Indag. Math.,} {\bf 25}, 634-653 (1963).

[45] J. Bernard and G. Letac, Construction d'evenements equiprobables et coefficients multinomiaux modulo $p^n$. {\it Illinois J. Math.,} {\bf 17}, N. 2, 317-332 (1997).

[46] H. Heyer, {\it Probability measures on locally compact groups.} Springer-Verlag (1977).

[47] K. Te\"oken, Limit theorems for probability measures on totally disconnected groups. {\it Semigroup Forum,} {\bf 58}, 69-84 (1999).

[48] R. Shah, Infinitely divisible measures on $p$-adic groups. {\it J. Theor. Prob.,} {\bf 4}, N. 2, 391-405 (1991).

[49] W. Hazod, Some new limit theorems for vector space- and group valued random variables. {\it J. Math. Soc.,} {\bf 93}, N. 4 (1999).

[50] W. Hazod and R. Shah, Semistable selfdecomposable laws on groups. {\it J. of Appl. Analysis,} {\bf 7}, N. 1, 1-22 (2001).

[51] A. Yu. Khrennikov, Interpretation of probability and their
$p$-adic extensions.  {\it Theory of Probability and Appl.,} {\bf
46}, N. 2, 311-325 (2001).

Albeverio S., Cianci R., De Grande-De Kimpe N., Khrennikov A.Yu.,
$p$-adic probability and an interpretation of negative probabilities
in quantum mechanics. {\it Russian J. of Math. Phys.,} {\bf 6}, 3-19
(1999).

[52] V. Maximov, Abstract models  of probability. Proc. of Int.
Conf. {\it Foundations of Probability and Physics,} ed. A.
Khrennikov, World Scientific Publishers (2001).

[53] A. Bendikov, L. Saloff-Coste, Brownian motion on compact groups in infinite dimension. Preprint of Cornell Univ. (2001).

[54] M. S. Barnett, Negative probability. {\it Proc. Cambridge Phil. Soc.,} {\bf 41}, 71-73 (1944).

[55] V. Yu. Krylov, Integration of analytic functionals with respect to signed measures. {\it Dokl. Akad. Nauk USSR,} {\bf 163}, 289-292 (1965).

[56] O. G. Smolyanov and A. Yu. Khrennikov, Central limit theorem for generalized measures on an infinite dimensional space. {\it Dokl. Akad. Nauk USSR,} {\bf 292}, N. 6, 1310-1313 (1987).

[57] A. Yu. Khrennikov, Integration with respect to generalized measures on linear topological spaces. {\it Trudy Mosovskogo Matematichaskogo Obshshestva,}
 {\bf 49}, 113-129 (1986) (English Translation in {\it Proc. Moscow Math. Soc.}).

 [58] A. Yu. Khrennikov, {\it Superanalysis.} Nauka, Moscow, (1997), in Russian. English translation: Kluwer Academic Publishers (1999).

 [59] D. Neuenschwander, {\it Comtes Rendus Acad. Sci. Paris,} Serie I, {\bf 330}, 1025-1030 (2000).

\end{document}